\documentclass[10pt]{article}

\usepackage[utf8]{inputenc}
\usepackage[letterpaper, margin = 1in]{geometry}

\usepackage{titlesec}
\titleformat{\section}[hang]{\large\bfseries}{\thesection}{0.5cm}{}
\titleformat{\subsection}[hang]{\normalsize\bfseries}{\thesubsection}{0.3cm}{}

\usepackage{graphicx}
\graphicspath{{Figures/}}
\usepackage{subcaption}

\usepackage{mathtools, amssymb, bm, amsthm}

\newtheorem{theorem}{Theorem}
\theoremstyle{definition}
\newtheorem{problem}{Problem}
\newtheorem{definition}{Definition}
\newtheorem{lemma}{Lemma}
\newtheorem{proposition}{Proposition}

\usepackage{array}
\usepackage{enumerate}
\usepackage{enumitem}
\usepackage{cite}
\usepackage{xcolor}
\usepackage{multirow}
\usepackage{booktabs}
\usepackage{algorithm}
\usepackage{algpseudocode}

\usepackage{hyperref}
\hypersetup{colorlinks = true, linkcolor = blue, citecolor = blue, urlcolor = blue}

\usepackage{tikz}
\usetikzlibrary{shapes, arrows, positioning, calc, intersections}
\tikzstyle{block} = [draw, thick, rectangle, fill=white!20, minimum height=3em, minimum width=6em]
\tikzstyle{square} = [draw, thick, rectangle, fill=white!20, minimum height = 3em, minimum width = 3em]
\tikzstyle{dotbox} = [draw, dashed, thick, rectangle, fill=white!20, minimum height = 3em, minimum width = 3em]
\tikzstyle{coord} = [coordinate]
    
\allowdisplaybreaks

\renewcommand\qedsymbol{$\blacksquare$}
\newcommand{\define}{ \stackrel{\Delta}{=} }
\newcommand{\bsmtx}{\left[ \begin{smallmatrix}}
\newcommand{\esmtx}{\end{smallmatrix} \right]} 

\title{\Large{\textbf{Analysis of Gradient Descent with Varying Step Sizes using\\Integral Quadratic Constraints}}}

\author{Ram Padmanabhan\thanks{Department of Electrical and Computer Engineering and Coordinated Science Laboratory, University of Illinois Urbana-Champaign, Urbana, IL 61801, USA. Email: \texttt{ramp3@illinois.edu}}~ and Peter Seiler\thanks{Department of Electrical Engineering and Computer Science, University of Michigan, Ann Arbor, MI 48109, USA. Email: \texttt{pseiler@umich.edu}}}
\date{}

\begin{document}

\maketitle

\begin{abstract}
The framework of Integral Quadratic Constraints (IQCs) is used to perform an analysis of gradient descent with varying step sizes. Two performance metrics are considered: convergence rate and noise amplification.  We assume that the step size is produced from a line search and varies in a  known interval. Modeling the algorithm as a linear, parameter-varying (LPV) system, we construct a parameterized linear matrix inequality (LMI) condition that certifies algorithm performance, which is solved using a result for polytopic LPV systems. Our results provide convergence rate guarantees when the step size lies within a restricted interval. Moreover, we recover existing rate bounds when this interval reduces to a single point, i.e. a constant step size. Finally, we note that the convergence rate depends only on the condition number of the problem. In contrast, the noise amplification performance depends on the individual values of the {strong convexity} and smoothness parameters, and varies inversely with them for a fixed condition number.
\end{abstract}

\section{Introduction} \label{sec:Introduction}
Convex optimization problems can be solved using numerous iterative algorithms, including gradient descent \cite{BV2004, NW2006}, accelerated \cite{BTP64, YN2004} and proximal gradient methods \cite{AB17, PB14}. The use of control theoretic methods to analyze such algorithms has a rich history. Our focus is on using the notion of Integral Quadratic Constraints (IQCs) from robust control theory. IQCs were first introduced by Yakubovich \cite{VAY92} in the context of imposing quadratic constraints on an infinite-horizon control problem in a Lur'e system, and imposing multiple constraints via the S-procedure \cite{VAY71}. 
Megretski and Rantzer \cite{MR97} unified the approach to robustness analysis of such systems using IQCs. 
Using the Kalman-Yakubovich-Popov (KYP) lemma \cite{AR96}, it was shown that verifying stability reduces to a linear matrix inequality (LMI) condition. Additional details can be found in \cite{Veenman16, Scherer22} with discussions connecting time- and frequency-domain IQCs \cite{PS2015} and discrete-time IQCs \cite{Hu17, BLR15, LL2016}. 

In \cite{DT13}, Drori and Teboulle developed a semidefinite program (SDP) approach to certify tight bounds for first-order optimization algorithms on strongly convex problems. Subsequent work by Lessard \emph{et al.} \cite{LL2016} developed a unified framework using IQCs for the analysis of first-order algorithms on strongly convex functions. Importantly, while the approach in \cite{DT13} scaled with the problem dimension, the use of IQCs circumvented this issue in \cite{LL2016}. A class of $\rho$-hard IQCs was introduced for characterizing convergence rates. Numerous IQCs for the gradient of strongly convex functions were presented, from simple sector bounds \cite{YN2004} to dynamic constraints using Zames-Falb IQCs \cite{HW2005}. Using these IQCs, certifying the convergence of first-order algorithms reduces to solving a small SDP independent of problem dimension. This framework has led to a significant body of subsequent work using IQCs to analyze and design optimization algorithms. This includes the analysis of the heavy-ball method \cite{BS2019}, the biased stochastic gradient method \cite{HSL21}, transient behavior of Nesterov's accelerated method \cite{MSJ23}, analysis of non-strongly convex problems \cite{Fazy18} and algorithm design \cite{SFL18, CHSL18, LS20}. A set of related case studies is available in \cite{LL22}.


Another important aspect of analyzing optimization algorithms is their robustness to noise, either in the iterate or in the gradient. Mohammadi \emph{et al.} \cite{MRJ19, MRJ21} examine the variance in the iterate error when iterates are perturbed by additive white noise, for both gradient descent and Nesterov's accelerated method. Based on this, a set of tradeoffs between noise amplification and convergence are derived in \cite{MRJ22} for strongly convex quadratic problems.

It must be noted that almost all prior work in this area assumes constant step sizes in the optimization algorithms. However, the use of constant step sizes usually requires additional knowledge of the problem, such as the Lipschitz constant $L$ of the gradient. Varying step sizes, often based on a line search \cite{LA66, PW69} are commonly implemented in optimization algorithms. 
While these algorithms are computationally more intensive, they do not require any prior knowledge about function parameters and do not result in a reduction in performance. In this article, we focus on the analysis of gradient descent with a varying step size based on the framework developed in \cite{LL2016}. We assume a line search has been carried out, and produces a step size that lies in a restricted interval about $\alpha = \frac{1}{L}$, a popular choice for gradient descent on strongly convex functions. While constant step size algorithms can be represented as linear, time-invariant (LTI) systems, a varying step size indicates that the algorithm is now a linear, parameter-varying (LPV) system, where the step size is a parameter affecting the dynamics. IQCs have been used for the analysis of LPV systems, including obtaining their worst-case $\mathcal{L}_2$ gain \cite{PF2015a, PF2016} and developing a robust synthesis algorithm \cite{WPF16}.

In our work, we combine the framework developed in \cite{LL2016} with the LPV analysis approaches discussed above, for analyzing gradient descent with a varying step size. {Using a result for polytopic LPV systems \cite{WYPB},} we obtain convergence rates for the parameter-varying algorithm. In addition to the analysis of convergence rate, we use the noise amplification framework developed in \cite{MRJ21} for constant step sizes, for our parameter-varying algorithm affected by gradient noise.

The remainder of this article is organized as follows. In Section \ref{sec:Preliminaries}, we present the two problems we examine, on convergence rate and noise amplification{, as well as some background on IQCs.} 
Our approach is described in Section \ref{sec:Varying} based on the characterization of the algorithm as an LPV system. This includes our main theorems, which extend prior constant step size results \cite{LL2016, MRJ19, MRJ21}. We also discuss how the implementation of the results can be simplified by tweaking the existing LMIs. Section \ref{sec:Results} presents our main results. First, we show that our approach recovers known constant step size expressions when the step size interval reduces to a single point. Further, we derive analytical expressions that characterize convergence and noise amplification for varying step size gradient descent based on our formulation. Finally, numerical results that are obtained by solving the LMIs are discussed. Section \ref{sec:Conclusion} provides concluding remarks.

\section{Preliminaries} \label{sec:Preliminaries}

\subsection{Notation} \label{sec:Notation}
Throughout this article, $I_n$ denotes the $n\times n$ identity matrix and $0_n$ denotes the $n\times n$ matrix of zeros. $\mathbb{Z}_+$ denotes the set of all non-negative integers, and $\mathbb{R}^n$ denotes the vector space of all $n$-tuples of real numbers. For a vector $x = [x_1, \ldots, x_n]^T \in \mathbb{R}^n$, $\|x\|$ denotes the Euclidean norm. $A\otimes B$ denotes the Kronecker product of two matrices $A$ and $B$. The Kronecker Delta function, denoted $\delta_j$ is defined as $\delta_j = 1$ if $j = 0$, and $\delta_j = 0$ otherwise.

\subsection{Problem Formulation} \label{sec:ProbForm}
Consider the following unconstrained optimization problem:
\begin{equation} \label{eq:unconopt}
	\min_{x\in \mathbb{R}^n} f(x),
\end{equation}
where $f:\mathbb{R}^n\to\mathbb{R}$. Let $\nabla f(x)$ denote the gradient of $f$ at $x$. We assume that $f$ is $m$-strongly convex and $\nabla f$ is $L$-Lipschitz continuous, i.e., for all $x, y \in \mathbb{R}^n$,
\begin{align}
	f(y) \geq &f(x) + \nabla f(x)^T(y-x) + \frac{m}{2}\|y-x\|^2 \label{eq:mcvx} \\
	\text{and }\hspace{0.5cm}&\left\|\nabla f(y) - \nabla f(x)\right\| \leq L\|y - x\|. \label{eq:Lipschitz}
\end{align}
The condition number of the function $f$ is $\kappa = \frac{L}{m}$. The class of functions $f$ that satisfy \eqref{eq:mcvx} and \eqref{eq:Lipschitz} is denoted $\mathcal{S}(m, L)$. If $f \in \mathcal{S}(m, L)$ there exists a unique solution $x^*$ to \eqref{eq:unconopt} \cite[Chapter 4]{BV2004}. We assume, without loss of generality by a coordinate shift, that this optimal point occurs at the origin, i.e. $x^* = 0$. 

In this article, we consider the use of gradient descent with a time-varying step size $\alpha_k$ as an iterative approach to solve the problem \eqref{eq:unconopt}. Initialized at some $x_0 \in \mathbb{R}^n$, gradient descent is characterized by the following update rule:
\begin{equation} \label{eq:gd-noisy}
    x_{k+1} = x_k - \alpha_k \left(\nabla f(x_k) + w_k\right).
\end{equation}
Here, $w_k$ is a perturbation used to model a noisy gradient. For the noise-free gradient method, $w_k = 0$ for all $k \in \mathbb{Z}_+$.

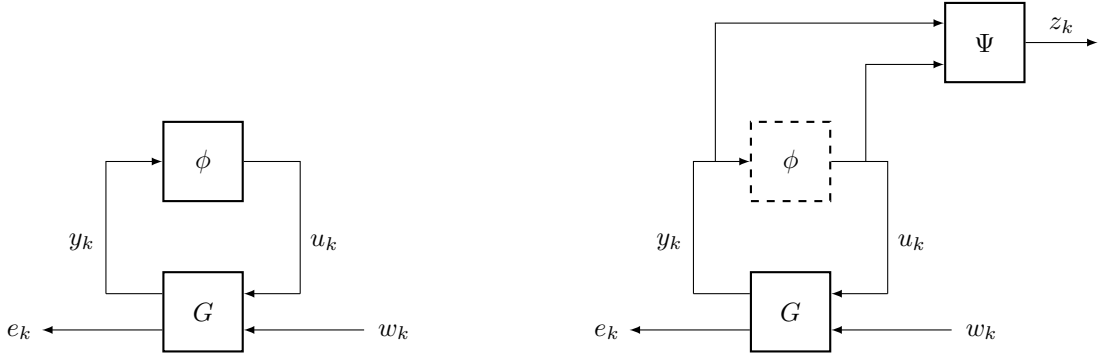
\begin{figure*}[!t]
\begin{subfigure}{0.49\textwidth}
\begin{center}
\begin{tikzpicture}[auto,>=latex',connect with angle/.style=
{to path={let \p1=(\tikztostart), 
\p2=(\tikztotarget) in -- ++({\x2-\x1-(\y1-\y2)*tan(#1-90)},0) -- (\tikztotarget)}}]

\node [square] (plant) {$G$};
\node [square, above of = plant, node distance = 2cm] (gradient) {$\phi$};
\node [coord, above left = -0.3cm and 0cm of plant] (n3) {};
\node [coord, above right = -0.3cm and 0cm of plant] (n4) {};
\node [coord, above left = -0.3cm and 0.75cm of plant] (n1) {};
\node [coord, above right = -0.3cm and 0.75cm of plant] (n2) {};
\node [coord, below left = -0.3cm and 0cm of plant] (n7) {};
\node [coord, below right = -0.3cm and 0cm of plant] (n8) {};
\node [coord, below left = -0.3cm and 1.6cm of plant] (n5) {};
\node [coord, below right = -0.3cm and 1.6cm of plant] (n6) {};

\draw[-latex] (n2) -- node{} (n4);
\draw[-] (n2) |- node[xshift = 0.32cm, yshift = -0.85cm]{$u_k$} (gradient);
\draw[-latex] (n1) |- node[xshift = -0.32cm, yshift = -1.3cm]{$y_k$} (gradient);
\draw[-] (n1) -- node{} (n3);
\draw[-latex] (n6) -- node[xshift = 1.2cm, yshift = 0.2cm]{$w_k$} (n8);
\draw[-latex] (n7) -- node[xshift = -1.1cm, yshift = 0.2cm]{$e_k$} (n5);

\end{tikzpicture}
\end{center}
\caption{A linear system $G$ in feedback with a nonlinear component $\phi$, affected by white stochastic noise $w_k$. In gradient descent, $\phi = \nabla f$.}
\end{subfigure}
\hfill
\begin{subfigure}{0.49\textwidth}
\begin{center}
\begin{tikzpicture}[auto,>=latex',connect with angle/.style=
{to path={let \p1=(\tikztostart), 
\p2=(\tikztotarget) in -- ++({\x2-\x1-(\y1-\y2)*tan(#1-90)},0) -- (\tikztotarget)}}]

\node [square] (plant) {$G$};
\node [dotbox, above of = plant, node distance = 2cm] (gradient) {$\phi$};
\node [square, above right = 2.5cm and 1.5cm of plant] (Psi) {$\Psi$};
\node [coord, above left = -0.3cm and 0cm of plant] (n3) {};
\node [coord, above right = -0.3cm and 0cm of plant] (n4) {};
\node [coord, above left = -0.3cm and 0.75cm of plant] (n1) {};
\node [coord, above right = -0.3cm and 0.75cm of plant] (n2) {};
\node [coord, below left = -0.3cm and 0cm of plant] (n7) {};
\node [coord, below right = -0.3cm and 0cm of plant] (n8) {};
\node [coord, below left = -0.3cm and 1.6cm of plant] (n5) {};
\node [coord, below right = -0.3cm and 1.6cm of plant] (n6) {};
\node [coord, left of = gradient, node distance = 1cm] (n31) {};
\node [coord, right of = gradient, node distance = 1cm] (n41) {};
\node [coord, right of = Psi, node distance = 1.5cm] (n51) {};
\node [coord, above right = 1.3cm and 1.5 cm of gradient] (n61) {};
\node [coord, above right = 0.75cm and 1.5 cm of gradient] (n71) {};

\draw[-latex] (n2) -- node{} (n4);
\draw[-] (n2) |- node[xshift = 0.32cm, yshift = -0.85cm]{$u_k$} (gradient);
\draw[-latex] (n1) |- node[xshift = -0.32cm, yshift = -1.3cm]{$y_k$} (gradient);
\draw[-] (n1) -- node{} (n3);
\draw[-latex] (n6) -- node[xshift = 1.2cm, yshift = 0.2cm]{$w_k$} (n8);
\draw[-latex] (n7) -- node[xshift = -1.1cm, yshift = 0.2cm]{$e_k$} (n5);
\draw[-latex] (n31) |- node{} (n61);
\draw[-latex] (n41) |- node{} (n71);
\draw[-latex] (Psi) -- node{$z_k$} (n51);

\end{tikzpicture}
\end{center}
\caption{Representation of the IQC framework. $z_k = \Psi(y_k, u_k)$, where $\Psi$ is a stable LTI system in the most general case.}
\end{subfigure}
\caption{$\phi$ is the nonlinear component we wish to analyze, and is replaced by the constraints it imposes on the input-output pair $(u, y)$. These are written as constraints on $z_k$.}
\label{fig:IQC}
\vspace{-1em}
\end{figure*}

We examine two fundamental performance metrics associated with the use of gradient descent:  (i) convergence rate and (ii) noise amplification. These metrics are defined as follows:

\begin{definition}[Convergence Rate] \label{def:CR}
The noise-free gradient method converges with a rate $\rho \in (0, 1)$ if there exists a constant {$\beta > 0$} such that:
\begin{equation} \label{eq:convergence}
	\|x_k\| \leq {\beta}\rho^k\|x_0\|
\end{equation}
for all $x_0 \in \mathbb{R}^n$ and for all $k \in \mathbb{Z}_+$.
\end{definition}

\begin{definition}[Noise Amplification] \label{def:NA}
Assume $w_k$ is additive white stochastic noise with zero mean and identity covariance matrix, i.e. $\mathbb{E}[w_k] = 0$ and $\mathbb{E}[w_kw_{l}^{T}] = I\delta_{k-l}$. The noise amplification of the algorithm \eqref{eq:gd-noisy} is characterized by the metric $\gamma$, where:
\begin{equation} \label{eq:metric}
	\gamma \define \left(\lim_{N\to\infty} \frac{1}{N} \sum_{k = 0}^{N}\mathbb{E}\left[\|x_k\|^2\right]\right)^{\frac{1}{2}} = \left(\lim_{N \to \infty} \mathbb{E}\left[\|x_N\|^2\right]\right)^{\frac{1}{2}}.
\end{equation}
$\gamma^2$ is the steady-state variance of the iterate $x_N$, which is also the steady-state variance of the iterate error $x_N - x^*$ as the optimal solution is $x^* = 0$. The limit above exists when the update \eqref{eq:gd-noisy} is stable.
\end{definition}

Our objective is to use the framework of Integral Quadratic Constraints (IQCs) to analyze these two performance metrics for gradient descent with a varying step size. We assume that $\alpha_k$ is a varying step size produced by some line search algorithm. Then, the two analysis problems we consider can be stated as follows:
\begin{problem}[Convergence Rate] \label{problem:CR}
For the noise-free gradient method, find the smallest possible convergence rate $\rho$ such that \eqref{eq:convergence} is satisfied.
\end{problem}
\begin{problem}[Noise Amplification] \label{problem:NA}
For the gradient method affected by gradient noise \eqref{eq:gd-noisy}, find the smallest possible bound on the steady-state iterate variance $\gamma^2$, and subsequently the metric $\gamma$ defined in \eqref{eq:metric}.
\end{problem}

\subsection{Integral Quadratic Constraints} \label{sec:IQCs}
The use of Integral Quadratic Constraints (IQCs) is motivated by the fact that eq. \eqref{eq:gd-noisy} can be written as a linear update rule separated from the gradient, which is a nonlinear component:
\begin{align} \label{eq:gd-noisy-2}
\begin{split}
    x_{k+1} &= x_k - \alpha_k u_k + \alpha_k w_k, \\
    y_k &= x_k, \\
    u_k &= \nabla f(y_k).
\end{split}
\end{align}
A general version of this configuration is shown in Fig.~\ref{fig:IQC}(a), where a linear system $G$ is in feedback with a static, memoryless nonlinear function $\phi$ such that $u = \phi(y)$. The dynamics of the linear system can, in general, be represented as:
{\begin{align}
\begin{split} \label{eq:LTI-G-noise}
    x_{k+1} &= A(\alpha_k) x_k + B_u(\alpha_k)u_k + B_w(\alpha_k)w_k, \\
    \begin{bmatrix} y_k \\ e_k \end{bmatrix} &= \begin{bmatrix} C_y(\alpha_k) \\ C_e(\alpha_k) \end{bmatrix}x_k.
\end{split}
\end{align}}
In this configuration, $e_k$ is a performance output. {While $G$ is linear, it may be time-varying or parameter-varying in which case the matrices $A$, $B_u$, $B_w$, $C_y$ and $C_e$ depend on the time index $k$ or a parameter $\alpha_k$.}

Analyzing the general interconnection in Fig.~\ref{fig:IQC}(a) is not straightforward due to the nonlinearity $\phi$. The IQC framework provides a convenient method to analyze this interconnection by replacing $\phi$ with the (usually quadratic) constraints it imposes on the input-output pair $(u, y)$. In Fig.~\ref{fig:IQC}(b), the general representation used for this framework is illustrated. $\Psi$ is a stable, LTI system that generates an auxiliary sequence $z_k$ from the sequences $y_k$ and $u_k$, i.e. $z = \Psi(y, u)$. Throughout this article, we consider the simplest case where $\Psi = \begin{bmatrix} 1 & 0 \\ 0 & 1 \end{bmatrix} \otimes I_n$, i.e. $z_k = \begin{bmatrix} y_k \\ u_k\end{bmatrix}$. $\Psi$ is simply a static map in this case. We now define two classes of IQCs that are useful in analyzing optimization algorithms.

\begin{definition} \label{def:Pointwise}
Consider the sequences $y_k$ and $u_k$, and let $z_k = \begin{bmatrix} y_k \\ u_k\end{bmatrix}$. The nonlinear function $u_k = \phi(y_k)$ satisfies the pointwise IQC defined by $M$ if for all $k \geq 0$,
\begin{equation} \label{eq:Pointwise}
z_{k}^{T}Mz_k \geq 0.
\end{equation}
\end{definition}

\begin{definition} \label{def:Rho}
Let $\rho > 0$ be given. Consider the sequences $y_k$ and $u_k$, and let $z_k = \begin{bmatrix} y_k \\ u_k\end{bmatrix}$. The nonlinear function $u_k = \phi(y_k)$ satisfies the $\rho$-hard IQC defined by $(M, \rho)$ if for all $T \geq 0$,
\begin{equation} \label{eq:Rho}
\sum_{k = 0}^{T}\rho^{-2k}z_{k}^{T}Mz_k \geq 0.
\end{equation}
\end{definition}

Note that if $\phi$ satisfies a pointwise IQC defined by $M$, it also satisfies the $\rho$-hard IQC defined by $(M, \rho)$ for all $\rho < 1$. This fact is useful in deriving convergence rates. 
If $z \in \mathbb{R}^{n_z}$, then $M \in \mathbb{R}^{n_z\times n_z}$, and is symmetric. Typically, $M$ is indefinite. Finally, more general classes of IQCs where $\Psi$ is a dynamical system are defined in \cite{LL2016}, which fall under the class of discrete-time Zames-Falb IQCs \cite{HW2005}.


For the case of gradient descent, the nonlinearity $\phi = \nabla f$ satisfies the following IQC (adapted from \cite{LL2016}), and this is used throughout this article.

\begin{lemma}[Sector IQC, \cite{LL2016}] \label{lem:Sector}
Let $f \in \mathcal{S}(m, L)$ and $\phi = \nabla f$. Then, $\nabla f$ satisfies the pointwise IQC defined by:
\begin{subequations} \label{eq:Sector}
\begin{equation}
	M = \begin{bmatrix} -2mL & (L+m) \\ (L+m) & -2 \end{bmatrix} \otimes I_n,
\end{equation}
\text{and the corresponding quadratic constraint is:}
\begin{equation}
    \begin{bmatrix} y_k \\ u_k \end{bmatrix}^T\begin{bmatrix} -2mLI_n & (L+m)I_n \\ (L+m)I_n & -2I_n \end{bmatrix}\begin{bmatrix} y_k \\ u_k \end{bmatrix} \geq 0.
\end{equation}
\end{subequations}
\end{lemma}

\section{Analysis with Varying Step Sizes} \label{sec:Varying}
In this section, we describe our approach to solve Problems \ref{problem:CR} and \ref{problem:NA} for gradient descent with a varying step size, using Integral Quadratic Constraints. This involves considering the algorithm \eqref{eq:gd-noisy} as a linear, parameter-varying (LPV) system. In Section \ref{sec:Varying-CR}, we present our main result on the convergence rate for the noise-free gradient method with a varying step size. In Section \ref{sec:Varying-NA}, we consider gradient descent affected by a noisy gradient and present our main result on noise amplification. Finally, Section \ref{sec:Implementation} provides details on numerical implementations for our results.

A line search is one of the most common methods leading to an optimization algorithm with varying step sizes \cite{NW2006, LA66, PW69}. In our setup, we assume that a line search has been carried out, and produces a step size at each time step $k$, denoted $\alpha_k$. Moreover, we assume that the step size satisfies the following condition:
\begin{equation} \label{eq:AlphaBounds}
	\underline{\alpha} = \frac{1}{cL} \leq \alpha_k \leq \frac{c}{L} = \overline{\alpha},
\end{equation}
where $c \geq 1$ is some constant characterizing the interval $\mathcal{A} = [\underline{\alpha}, \overline{\alpha}]$. In other words, the step size $\alpha$ is not constant in time, but takes a value in this interval at each time step $k$. This interval is also geometrically centered about the popular step size choice $\alpha = 1/L$. While conditions on the step size after carrying out a Wolfe line search can be derived, the resulting bounds can be too conservative to be useful, particularly for high condition numbers.

Under the assumption \eqref{eq:AlphaBounds}, the gradient algorithm \eqref{eq:gd-noisy} can be written as an LPV system:
\begin{align} \label{eq:LPV-G-noisy}
\begin{split}
    x_{k+1} &= Ax_k + B_u(\alpha_k)u_k + B_w(\alpha_k)w_k, \\
    \begin{bmatrix} y_k \\ e_k \end{bmatrix} &= \begin{bmatrix} C_y \\ C_e \end{bmatrix}x_k,
\end{split}
\end{align}
where
\begin{equation} \label{eq:gd-noisy-var-matrices}
    A = I_n; \hspace{0.15em} B_u(\alpha_k) = -\alpha_k I_n; \hspace{0.15em} B_w(\alpha_k) = \alpha_k I_n; \hspace{0.15em} C_y = C_e = I_n,
\end{equation}
and the matrices $B_u$ and $B_w$ depend on the parameter $\alpha_k \in \mathcal{A}$. 

\subsection{Convergence Rate} \label{sec:Varying-CR}
We first consider the noise-free case, where $w_k = 0$ for all $k \in \mathbb{Z}_+$. The following theorem provides a method to characterize convergence rates for gradient descent with a varying step size satisfying \eqref{eq:AlphaBounds}, and is an extension of the main result in \cite{LL2016}. 

\begin{theorem}[Problem \ref{problem:CR}, Convergence Rate] \label{thm:VaryingConv}
Suppose the nonlinear function $u = \phi(y)$ satisfies a pointwise IQC defined by $M$ as given in Definition \ref{def:Pointwise}. Suppose that there exists a positive definite matrix $P$, non-negative scalar $\lambda$, and scalar $\rho \in [0,1]$ such that the following LMI is feasible for all $\alpha \in \mathcal{A}$: 
\begin{align}
&\begin{bmatrix}
A^TPA - \rho^2P & A^TPB_u(\alpha) \\ B_{u}^{T}(\alpha)PA & B_{u}^{T}(\alpha)PB_u(\alpha)
\end{bmatrix} + \lambda \begin{bmatrix} C_y & 0_n \\ 0_n & I_n \end{bmatrix}^T M \begin{bmatrix} C_y & 0_n \\ 0_n & I_n \end{bmatrix} \preceq 0, \label{eq:VaryingConvLMI}
\end{align}
Then, for any $x_0$ we have:
\begin{equation} \label{eq:varying-expconv}
	\|x_k\| \leq \sqrt{\textup{cond}(P)} \rho^k \|x_0\|.
\end{equation}
\end{theorem}
\begin{proof}
{The proof closely follows the arguments in \cite[Theorem 4]{LL2016}, and is omitted to conserve space.}
\end{proof}

\subsection{Noise Amplification} \label{sec:Varying-NA}
Now consider the gradient algorithm affected by gradient noise, as written in \eqref{eq:LPV-G-noisy} and \eqref{eq:gd-noisy-var-matrices}. The following theorem characterizes the noise amplification metric $\gamma$ for this algorithm, and extends the constant step size result in \cite{MRJ19, MRJ21}.

\begin{theorem}[Problem \ref{problem:NA}, Noise Amplification] \label{thm:VaryingNA}
Suppose the nonlinear function $u = \phi(y)$ satisfies a pointwise IQC defined by $M$ as given in Definition \ref{def:Pointwise}. Suppose that there exists a positive definite matrix $P$ and a non-negative scalar $\lambda$ such that the following LMI is feasible for all $\alpha \in \mathcal{A}$: 
\begin{align} 
	&\begin{bmatrix} A^TPA - P + C_{e}^{T}C_e & A^TPB_u(\alpha) \\ B_{u}^{T}(\alpha)PA & B_{u}^{T}(\alpha)PB_u(\alpha) \end{bmatrix} + \lambda \begin{bmatrix} C_y & 0_n \\ 0_n & I_n \end{bmatrix}^TM\begin{bmatrix} C_y & 0_n \\ 0_n & I_n \end{bmatrix} \preceq 0. \label{eq:VaryingNALMI}
\end{align}
Then, the metric $\gamma$ is bounded by:
\begin{equation} \label{eq:gamma-1-varying}
    \gamma \leq \sup_{\alpha \in \mathcal{A}} \Big(\mathrm{tr}\left(B_{w}^{T}(\alpha)PB_w(\alpha)\right)\Big)^{\frac{1}{2}}.
\end{equation}
\end{theorem}
\begin{proof}
Define a storage function $V(x_k) = x_{k}^{T}Px_k.$ The LMI \eqref{eq:VaryingNALMI} must hold for $\alpha = \alpha_k$. Following similar steps to \cite[Lemma 1]{MRJ21}, we can show that:
\begin{align}
	\frac{1}{N} \sum_{k = 0}^{N} \mathbb{E}[\|e_k\|^2] \leq &\frac{1}{N} \mathbb{E}[V(x_0) - V(x_{N+1})] + \mathrm{tr}\left(B_{w}^{T}(\alpha_k)PB_w(\alpha_k)\right), \label{eq:Ee}
\end{align}
where $\mathbb{E}[.]$ is over different realizations of $w_k$. Note that
$$
\mathrm{tr}\left(B_{w}^{T}(\alpha_k)PB_w(\alpha_k)\right) \leq \sup_{\alpha \in \mathcal{A}} \mathrm{tr}\left(B_{w}^{T}(\alpha)PB_w(\alpha)\right).
$$
Using this in \eqref{eq:Ee}, taking the limit as $N \to \infty$, using the definition of $\gamma$ in \eqref{eq:metric} and noting that in our setup, the performance output $e_k$ is the state $x_k$, the result \eqref{eq:gamma-1-varying} follows.
\end{proof}

We note here that a common approach to reduce conservatism in parameter-varying LMI problems such as \eqref{eq:VaryingConvLMI} and \eqref{eq:VaryingNALMI} is to use a parameter-dependent matrix $P(\alpha)$ instead of a constant matrix $P$. This is used to address bounded rates of variation of the parameter, as discussed in \cite{WYPB}. However, in our setup we do not assume any bounds on the rate of variation of $\alpha$. Our only assumption is that $\alpha_k$ satisfies the constraint \eqref{eq:AlphaBounds} at each time step $k$, but can vary at any rate between these quantities. Thus, we consider only a constant matrix $P$.

\subsection{Numerical Implementation} \label{sec:Implementation}
We now discuss a few details on implementing the LMIs in \eqref{eq:VaryingConvLMI} and \eqref{eq:VaryingNALMI}. Each LMI is actually an infinite family of LMIs, as they must be satisfied for each $\alpha \in \mathcal{A}$. We can simplify this to a finite number of constraints using a result for parameterized LMIs \cite{WYPB}. Note how the LMIs \eqref{eq:VaryingConvLMI} and \eqref{eq:VaryingNALMI} can be written as:
\begin{equation} \label{eq:LMISimplified-1}
    \begin{bmatrix} A^T \\ B_{u}^{T}(\alpha) \end{bmatrix}P\begin{bmatrix} A & B_u(\alpha) \end{bmatrix} + X(P, \lambda) \preceq 0,
\end{equation}
for an appropriate matrix $X(P, \lambda)$ that is an affine function of the decision variables $P$ and $\lambda$. By Schur complements, this is equivalent to:
\begin{equation} \label{eq:LMISimplified-2}
    \left[\begin{array}{c|c}
        \multirow{2}{*}{$X(P, \lambda)$} & A^T \\ & B_{u}^{T}(\alpha) \\ \hline A~~~~B_u(\alpha) & -P^{-1}
    \end{array}\right] \preceq 0.
\end{equation}
Multiply the above equation on the left and right by the symmetric, block-diagonal matrix $\bsmtx I_{2n} & \bm{0} \\ \bm{0}^T & P \esmtx$ where $\bm{0}$ denotes the zero matrix of appropriate dimensions. Then, \eqref{eq:LMISimplified-2} is equivalent to:
\begin{equation} \label{eq:LMISimplified-3}
    \left[\begin{array}{c|c}
        \multirow{2}{*}{$X(P, \lambda)$} & A^TP \\ & B_{u}^{T}(\alpha)P \\ \hline PA~~~~PB_u(\alpha) & -P
    \end{array}\right] \preceq 0.
\end{equation}
Note that \eqref{eq:LMISimplified-3} is affine in $P$ and $\lambda$. Furthermore, in contrast to the form of the original LMIs \eqref{eq:VaryingConvLMI}, \eqref{eq:VaryingNALMI}, the parameter $\alpha$ now enters affinely in the above expression since $B_u(\alpha) = -\alpha I_n$. Thus, it is sufficient to check the constraints \eqref{eq:VaryingConvLMI}, \eqref{eq:VaryingNALMI} at the end points $\underline{\alpha}$ and $\overline{\alpha}$ \cite{WYPB, PF2015a}. In other words, there exist $P$ and $\lambda$ such that \eqref{eq:LMISimplified-3} holds for all $\alpha \in \mathcal{A}$ if and only if there exist $P$ and $\lambda$ such that \eqref{eq:LMISimplified-3} holds for $\alpha \in \{\underline{\alpha}, \overline{\alpha}\}$. This proves particularly useful in deriving analytical solutions of \eqref{eq:VaryingConvLMI} and \eqref{eq:VaryingNALMI}, {and is a standard technique for polytopic LPV systems such as the ones we consider in \eqref{eq:LPV-G-noisy} and \eqref{eq:gd-noisy-var-matrices}.}

Then, the convergence rate problem can be written as:
\begin{align}
&\min_{\rho, P \succ 0, \lambda \geq 0}~~ \rho^2 \nonumber \\
\text{s.t. }~ &\text{LMI \eqref{eq:VaryingConvLMI} holds for $\underline{\alpha}$ and $\overline{\alpha}$} \label{eq:CR-numerical}
\end{align}
This is a bilinear problem due to the presence of the term $\rho^2 P$ in the LMI \eqref{eq:VaryingConvLMI}. However, it is quasiconvex and is known as a generalized eigenvalue problem \cite{Boyd93}. While there exist special solvers for such problems, it can also be solved via bisection on $\rho^2$, checking the feasibility of the constraints \eqref{eq:CR-numerical} at each step. The noise amplification problem can be written as:
\begin{align}
&\min_{P \succ 0, \lambda \geq 0, \gamma > 0}~~ \gamma^2 \nonumber \\
\text{s.t. }~ &\text{LMI \eqref{eq:VaryingNALMI} holds for $\underline{\alpha}$ and $\overline{\alpha}$, } \nonumber \\
\text{and }~ &\gamma^2 \geq \mathrm{tr}\left(B_{w}^{T}(\alpha)PB_w(\alpha)\right) \text{ at } \underline{\alpha} \text{ and } \overline{\alpha} \label{eq:SDP-gamma}
\end{align}
This problem is an SDP and can be solved using freely available solvers.

\section{Results} \label{sec:Results}
In this section, we present the results of our study on gradient descent with varying step sizes. In Section \ref{sec:Const} we discuss the reduction to the constant step size case by setting the interval constant $c = 1$, and show that prior results on convergence rate and noise amplification are recovered. We also present some insights that are particular to the gradient noise setting, including a tradeoff between convergence rate and noise amplification for strongly convex functions. In Section \ref{sec:VaryingResults} we first present analytical expressions for the convergence rate and noise amplification metric as functions of the condition number $\kappa$ and the interval constant $c$, based on the parameter-varying LMIs \eqref{eq:VaryingConvLMI} and \eqref{eq:VaryingNALMI}. Finally, we present a set of numerical results based on Theorems \ref{thm:VaryingConv} and \ref{thm:VaryingNA}.

\subsection{Reduction to a Constant Step Size} \label{sec:Const}
When $c = 1$, the step size interval $\mathcal{A}$ reduces to a single point, resulting in a constant step size $\alpha$. The two LMIs \eqref{eq:VaryingConvLMI} and \eqref{eq:VaryingNALMI} reduce to a single condition for $\alpha = \underline{\alpha} = \overline{\alpha}$.

\begin{figure}[!t]
    \centering
    \includegraphics[width = 0.48\textwidth]{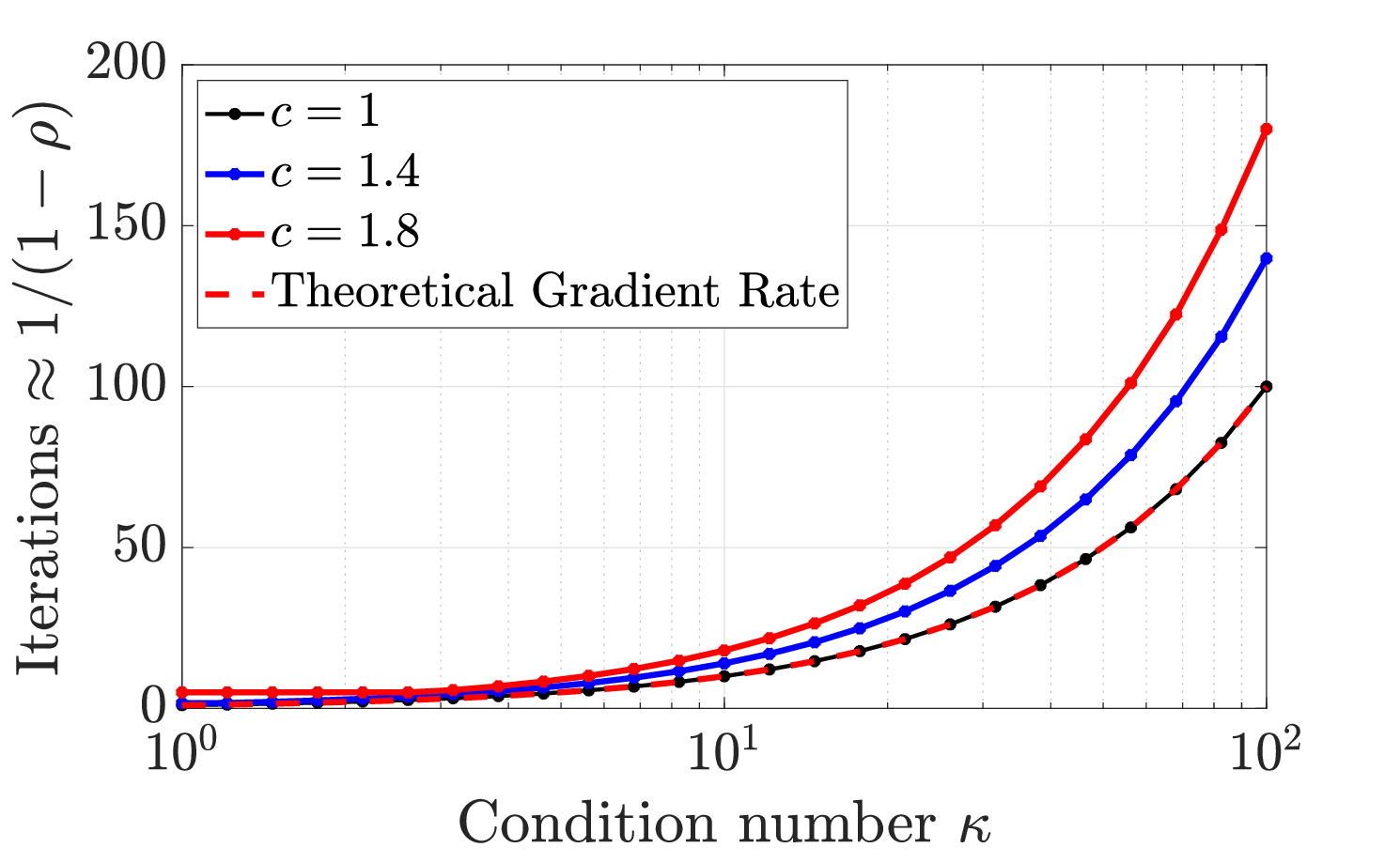}
    \caption{The approximate number of iterations to convergence as a function of the condition number for gradient descent with a varying step size characterized by $c$.}
    \label{fig:gd-varying-iter-conv}
    \vspace{-1em}
\end{figure}

For the convergence rate problem, substituting \eqref{eq:gd-noisy-var-matrices} and \eqref{eq:Sector} in \eqref{eq:VaryingConvLMI} when $c = 1$ produces the following:
\begin{equation} \label{eq:CCRSimplified}
\begin{bmatrix} (1-\rho^2) & -\alpha \\ -\alpha & \alpha^2 \end{bmatrix} + 
\lambda \begin{bmatrix} -2mL & L+m \\ L+m & -2 \end{bmatrix} \preceq 0,
\end{equation}
which is the same LMI obtained in \cite{LL2016}. This follows from a dimensionality reduction argument described therein. An analytical solution can be obtained using Schur complements:
\begin{equation} \label{eq:const-rho}
    \rho^* = \max\{|1-\alpha m|, |1-\alpha L|\} = 
    \begin{dcases}
        (1-\alpha m), \hspace{-1.5em} &\alpha \leq \frac{2}{L+m}, \\
        (\alpha L - 1), &\alpha > \frac{2}{L+m}
    \end{dcases},
\end{equation}
although constant step sizes larger than $\alpha = \frac{2}{L+m}$ are not common and do not provide an improvement in performance. For $\alpha = 1/L$, we note that the known convergence rate $\rho^* = 1-\frac{1}{\kappa}$ is recovered using \eqref{eq:const-rho}. {Figure~\ref{fig:gd-varying-iter-conv} illustrates the approximate number of iterations to convergence, given by $\frac{1}{(1-\rho)}$, as a function of the condition number $\kappa$ for different values of the interval constant $c$. We discuss the $c = 1$ case here, and other values of $c$ in Section \ref{sec:VaryingResults}. As we expect, the black curve for $c = 1$ coincides with the dashed red curve denoting the theoretical number of iterations $\frac{1}{(1-\rho^*)} = \kappa$ for the constant step size $\alpha = 1/L$. Setting $c = 1$ thus recovers the known constant step size result for the convergence of gradient descent, based on the family of LMIs in Theorem \ref{thm:VaryingConv}.}

For the noise amplification problem, the constant step size case was discussed in \cite{MRJ19, MRJ21} using the LMI \eqref{eq:VaryingNALMI} for a single point $\alpha$, under the iterate noise setting where $B_w = \sigma I_n$ and $\sigma$ was the noise magnitude. 
Here, we consider a noisy gradient for which the dynamics are characterized by \eqref{eq:gd-noisy-var-matrices} and where $B_w(\alpha) = \alpha I_n$ when $c = 1$. In Theorem \ref{thm:VaryingNA}, from \eqref{eq:gamma-1-varying}:
\begin{equation} \label{eq:gamma-1}
    \gamma \leq \big(\mathrm{tr}\left(B_{w}^{T}PB_w\right)\big)^{\frac{1}{2}},
\end{equation}
where $B_w$ is a constant for a given $\alpha$. Following identical steps to \cite{MRJ19}, the metric $\gamma$ satisfies the {following bound:}
\begin{equation} \label{eq:gamma-star}
    \gamma \leq \gamma^* = \alpha\sqrt{\frac{n}{(1-\rho^2)}},
\end{equation}
where $\rho$ is the convergence rate corresponding to the step size $\alpha$, obtained from \eqref{eq:const-rho} and $\gamma^*$ is the best upper bound on the noise amplification metric $\gamma$. For $\alpha = 1/L$, we have:
\begin{equation} \label{eq:gamma-popular}
    \gamma^* = \frac{1}{m}\sqrt{\frac{n}{2\kappa-1}} = \frac{1}{L}\sqrt{\frac{n}{2\kappa^{-1}-\kappa^{-2}}},
\end{equation}
and similarly if $\alpha = \frac{2}{L+m}$, we have:
\begin{equation} \label{eq:gamma-optimal}
    \gamma^* = \frac{1}{m}\sqrt{\frac{n}{\kappa}} = \frac{1}{L}\sqrt{\frac{n}{\kappa^{-1}}}.
\end{equation}
Notice how $\gamma^*$ depends separately on $m$ and $L$, and not just on the condition number $\kappa$. Further, note that this upper bound varies inversely with $m$ and $L$ for a fixed $\kappa$. Gradient noise amplification is thus an inherently different property from convergence rate, in that it depends individually on $m$ and $L$ and not just on $\kappa$. In particular, for a given condition number, $\gamma^*$ can take different values based on the values of $m$ and $L$.

\begin{figure*}[!t]
\centering
\begin{subfigure}{0.33\textwidth}
\centering
\includegraphics[width = \textwidth]{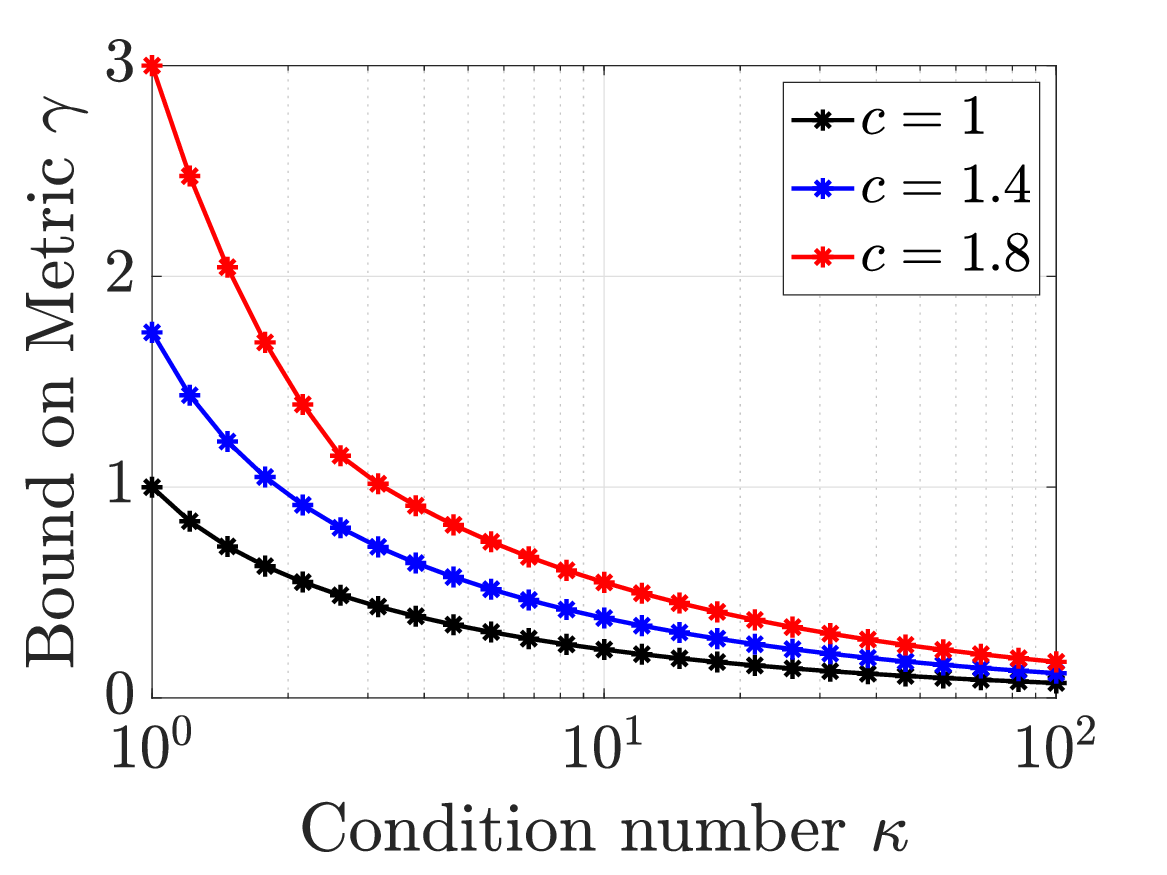}
\caption{Varying $L$ with $m = 1$.}
\end{subfigure}
\hspace{0.4cm}
\begin{subfigure}{0.33\textwidth}
\centering
\includegraphics[width = \textwidth]{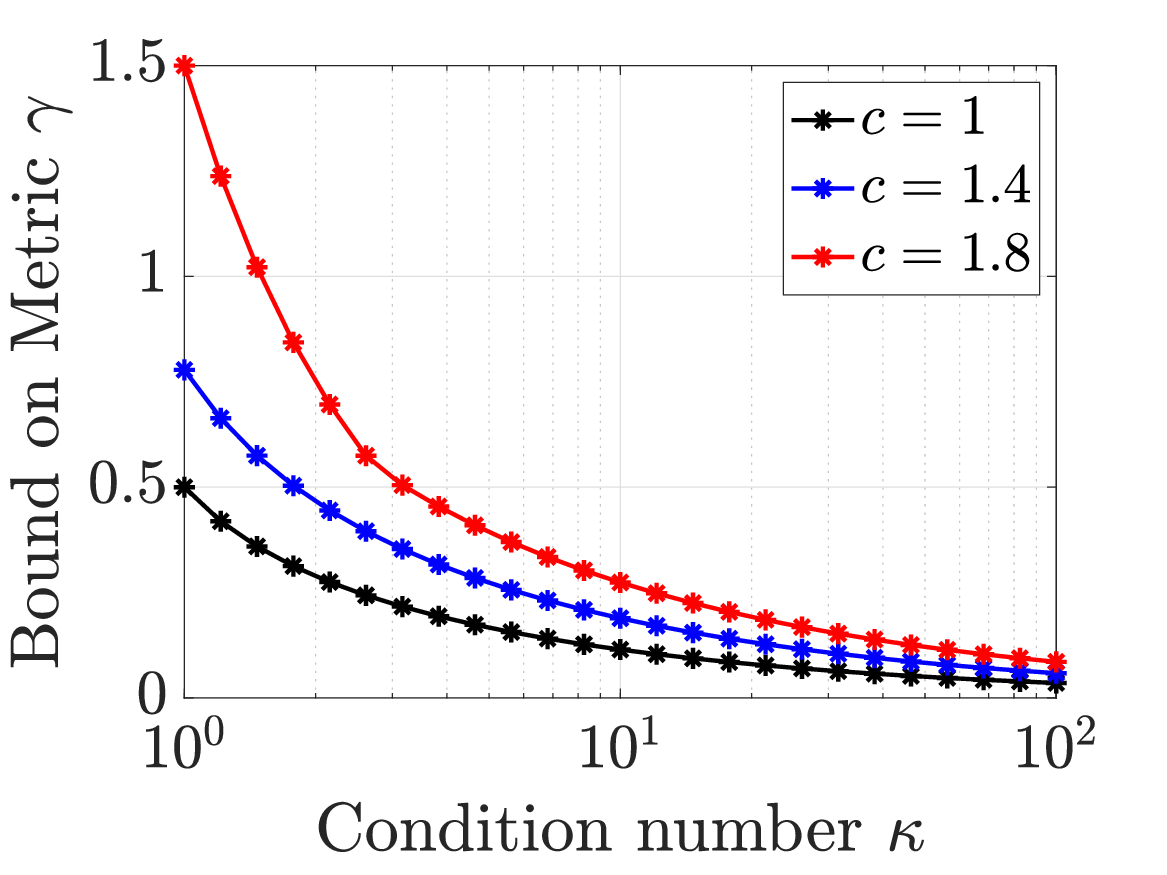}
\caption{Varying $L$ with $m = 2$.}
\end{subfigure}
\\
\begin{subfigure}{0.33\textwidth}
\centering
\includegraphics[width = \textwidth]{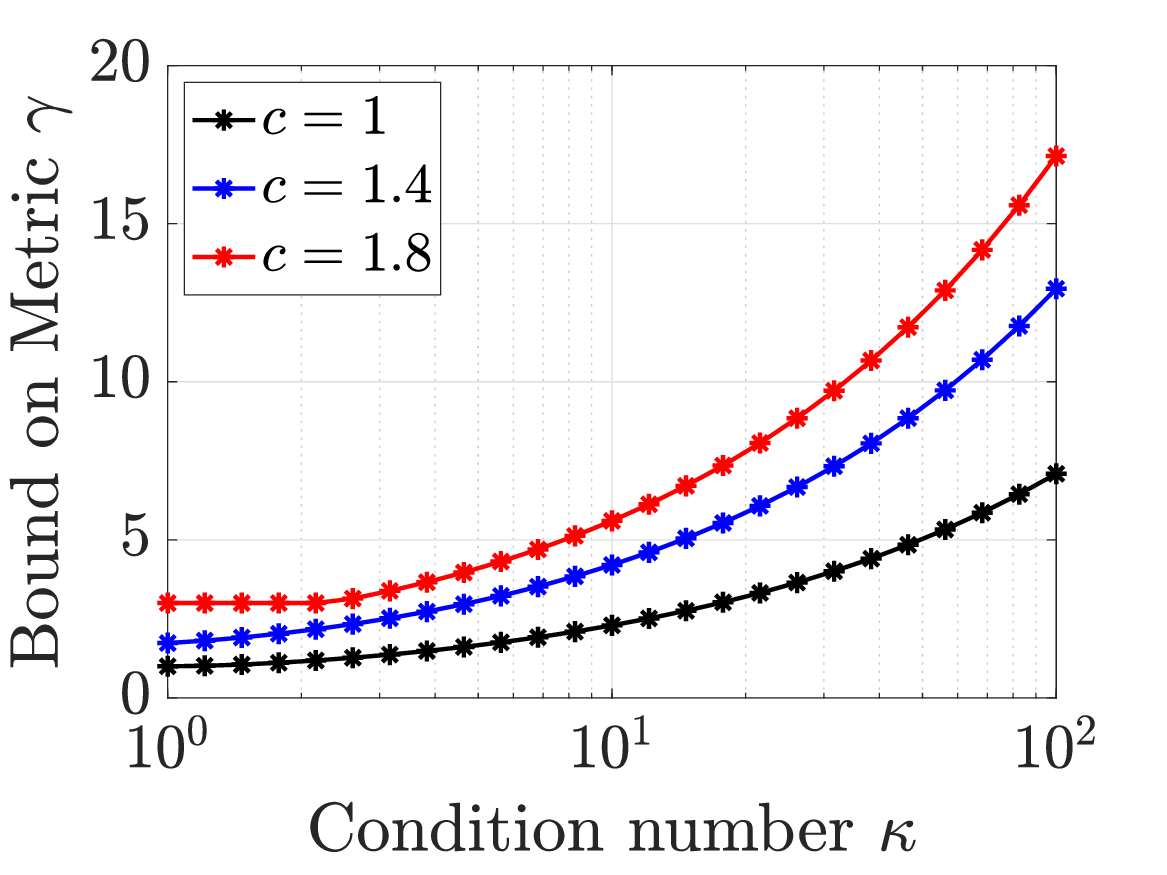}
\caption{Varying $m$ with $L = 1$.}
\end{subfigure}
\hspace{0.4cm}
\begin{subfigure}{0.33\textwidth}
\centering
\includegraphics[width = \textwidth]{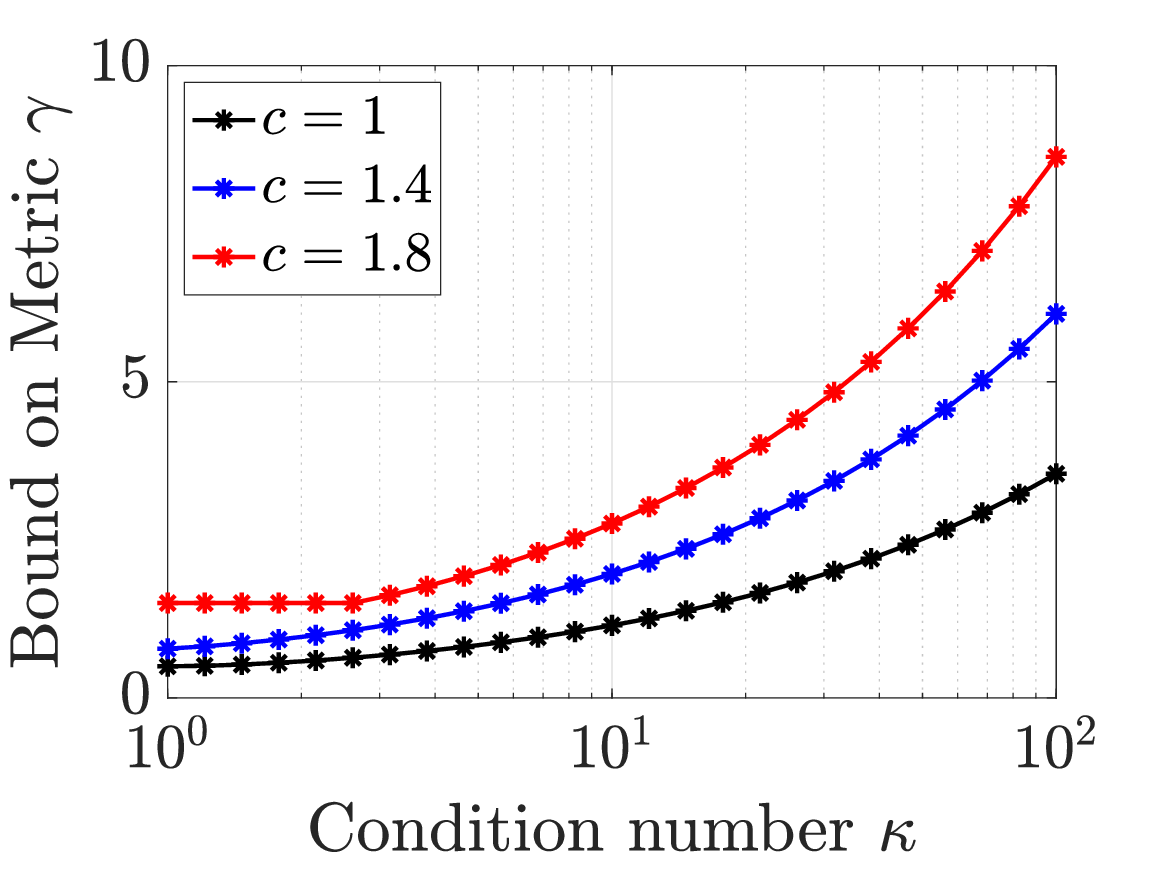}
\caption{Varying $m$ with $L = 2$.}
\end{subfigure}

\caption{The upper bound on the noise amplification metric as a function of the condition number for gradient descent with a varying step size characterized by $c$.}
\label{fig:gd-varying-na}
\vspace{-1.5em}
\end{figure*}

Finally, these expressions are very closely related to the corresponding expressions in \cite{MRJ19} for iterate noise amplification, and replacing $\alpha$ by $\sigma = 1$ in the above equations recovers the expressions obtained in \cite{MRJ19}. Thus, setting $c = 1$ can recover known expressions for the noise amplification metric for constant step sizes. However, in the iterate noise setting, the best upper bound on $\gamma$ no longer depends separately on $m$ and $L$, and depends only on the condition number $\kappa$, unlike the gradient noise setting discussed above.

\begin{figure}[!t]
    \centering
    \includegraphics[width = 0.48\textwidth]{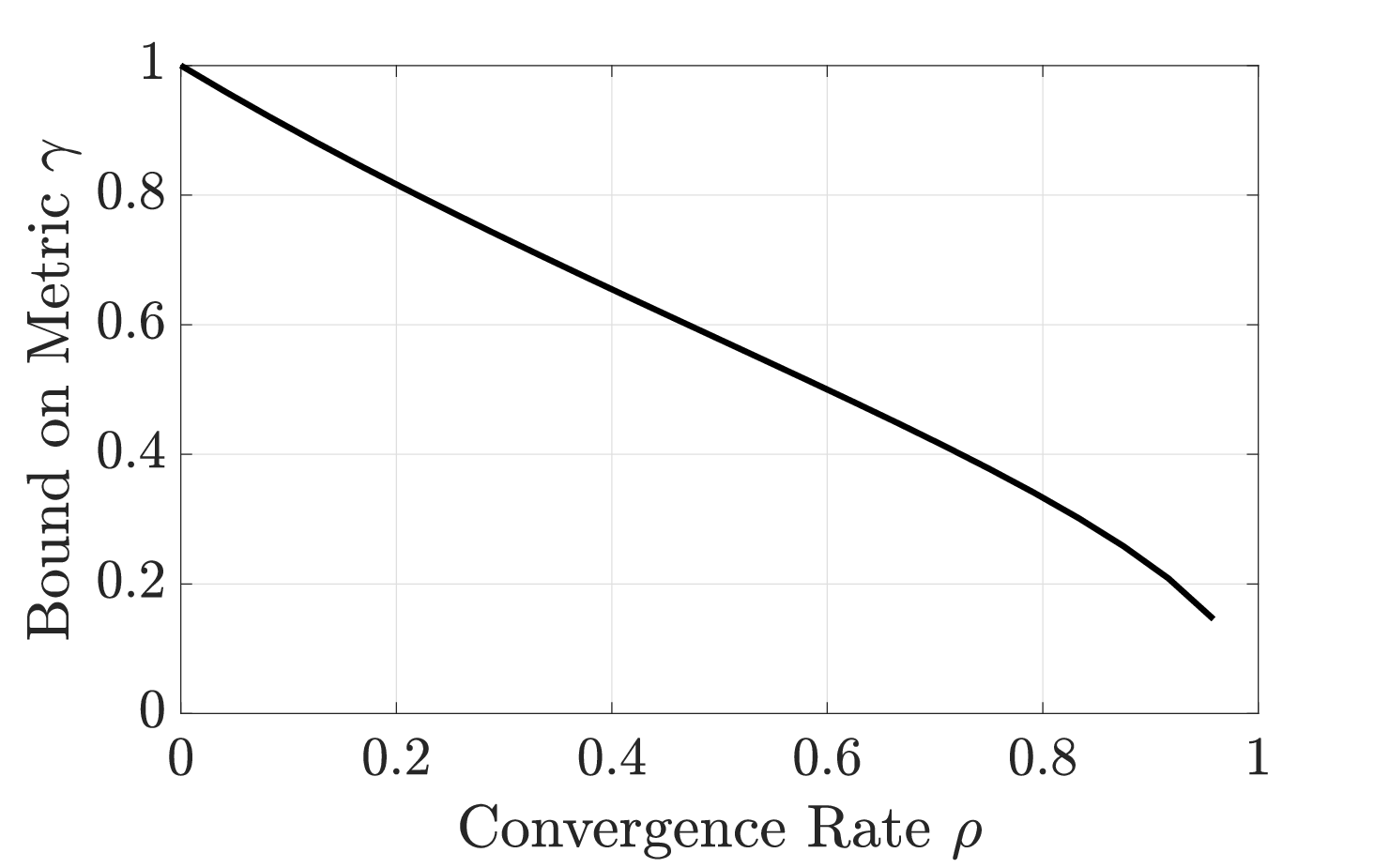}
    \caption{Tradeoff between noise amplification $\gamma$ and convergence rate $\rho$, based on \eqref{eq:tradeoff}. A `faster' algorithm has a larger value of metric $\gamma$, and is thus more sensitive to noise. In this figure, problem dimension $n = 1$ and strong convexity parameter $m = 1$.}
    \label{fig:tradeoff}
    \vspace{-1.5em}
\end{figure}

Figure~\ref{fig:gd-varying-na} illustrates the variation of $\gamma^*$ with condition number $\kappa$ for different values of the interval constant $c$, and fixing one of the two parameters $m$ and $L$. As before, we discuss the $c = 1$ case here, and other values of $c$ in Section \ref{sec:VaryingResults}. When $L$ is varied for particular values of $m$, the upper bound decreases as $L$ and $\kappa$ increase. This is shown in Figures~\ref{fig:gd-varying-na}(a) and (c). When $m$ is varied for particular values of $L$, $\gamma$ increases as $m$ decreases and $\kappa$ increases. This is shown in Figures~\ref{fig:gd-varying-na}(b) and (d). Furthermore, when $m$ or $L$ are made twice as large, the corresponding value of $\gamma^*$ is half as large for all condition numbers, as seen by comparing Figure~\ref{fig:gd-varying-na}(a) with \ref{fig:gd-varying-na}(b) and Figure~\ref{fig:gd-varying-na}(c) with \ref{fig:gd-varying-na}(d). This illustrates that $\gamma$ varies inversely with $m$ and $L$ for a given condition number.

We now discuss a tradeoff between convergence rate and noise amplification that arises from \eqref{eq:const-rho} and \eqref{eq:gamma-star}. {A cursory examination of \eqref{eq:gamma-star} seems to indicate that the noise amplification must increase as the convergence rate increases, i.e. worsens. However, note that $\gamma^*$ also depends on the step size $\alpha$, which is related to the convergence rate through \eqref{eq:const-rho}.} For $\alpha \leq \frac{2}{L+m}$, we have $\rho = (1-\alpha m)$ or $\alpha = \frac{1-\rho}{m}$. Then,
\begin{equation} \label{eq:tradeoff}
    \gamma^* = \frac{(1-\rho)}{m}\sqrt{\frac{n}{(1-\rho^2)}} = \frac{1}{m}\sqrt{\frac{n}{(1+\rho)}}.
\end{equation}
For $0 \leq \rho < 1$, $\gamma^*$ decreases with increasing $\rho$ as shown in Fig.~\ref{fig:tradeoff}. Thus, a choice of step size that results in faster convergence also is more sensitive to noise. This is closely related to a recent tradeoff result for accelerated momentum-based algorithms \cite{MRJ22} on strongly convex quadratic problems. However, our result here applies to the more general class of strongly convex functions.

\subsection{Results for a Varying Step Size} \label{sec:VaryingResults}
We now discuss our results for the convergence and noise amplification of gradient descent with a varying step size, based on Theorems \ref{thm:VaryingConv} and \ref{thm:VaryingNA}. First, we present analytical expressions for the rate of convergence and noise amplification metric for varying step sizes, and then discuss a set of numerical results by solving the problems \eqref{eq:CR-numerical} and \eqref{eq:SDP-gamma}.

The following propositions provide analytical expressions for the convergence rate and noise amplification metric with varying step sizes in gradient descent. {The proofs are provided in the appendix,} and rely on the fact that it is sufficient to check the families of constraints \eqref{eq:VaryingConvLMI} and \eqref{eq:VaryingNALMI} only at the end points $\underline{\alpha}$ and $\overline{\alpha}$ of the step size interval $\mathcal{A}$. 

\begin{proposition} \label{prop:CR}
The best upper bound on the convergence rate of varying step size gradient descent, based on the solution of the problem \eqref{eq:CR-numerical} for $1 \leq c < 2$ is given by:
\begin{equation} \label{eq:varying-rho}
    \rho_{\mathrm{varying}} = 
    \begin{dcases}
        c-1, \hspace{1em} &\kappa \leq \frac{1}{c(2-c)} \\
        1-\frac{1}{c\kappa}, &\kappa > \frac{1}{c(2-c)} 
    \end{dcases},
\end{equation}
depending on the interval constant $c$ and the condition number $\kappa$ as above. 
\end{proposition}

\begin{proposition} \label{prop:NA}
The noise amplification metric $\gamma_{\mathrm{varying}}$ for varying step size gradient descent, based on the solution of the SDP \eqref{eq:SDP-gamma} for $1 \leq c < 2$ satisfies:
\begin{equation} \label{eq:varying-gamma}
    \gamma_{\mathrm{varying}} \geq \gamma^{*}_{\mathrm{varying}} = 
    \begin{dcases}
        \frac{c}{L}\sqrt{\frac{n}{2c-c^2}}, \hspace{1em} &\kappa \leq \frac{c}{2-c} \\
        \frac{c}{m}\sqrt{\frac{n}{2c\kappa - c^2}}, &\kappa > \frac{c}{2-c}
    \end{dcases},
\end{equation}
depending on the interval constant $c$ and the condition number $\kappa$ as above. 
\end{proposition}

We now present a set of numerical results based on Theorems \ref{thm:VaryingConv} and \ref{thm:VaryingNA}, which characterize how the convergence rate $\rho$ and noise amplification metric $\gamma$ vary with the condition number $\kappa$ as well as the interval constant $c$. We choose three values of $c$ to test: $c = 1$, $c = 1.4$ and $c = 1.8$. Note that $c = 1$ implies that the interval $\mathcal{A}$ reduces to a point $\alpha_k \equiv \alpha = 1/L$, a constant step size, which was discussed in Section \ref{sec:Const}.

In the convergence rate results presented here, we use $m = 1$ throughout. However, the results are unchanged for other values of $m$ as the convergence rate depends only on the condition number $\kappa$ as seen in \eqref{eq:varying-rho}. {We now revisit Fig.~\ref{fig:gd-varying-iter-conv}, and discuss results for larger values of $c$. We first note that our results are consistent with the analytical expression derived in Proposition~\ref{prop:CR} for different values of $c$. As discussed earlier, setting $c = 1$ also recovers the theoretical number of iterations for convergence when $\alpha = 1/L$. For larger values of $c$, we observe that convergence is generally slower as the difference between $\underline{\alpha}$ and $\overline{\alpha}$ is larger, leading to a larger interval $\mathcal{A}$ and introducing some conservatism in the rate bound. Further note that if $c \geq 2$, convergence would not be guaranteed based on our approach to solve the LMI in \eqref{eq:VaryingConvLMI}.}


For the noise amplification results, we test different values of $m$ and $L$, fixing one of the two parameters. These results are shown in Fig.~\ref{fig:gd-varying-na}. As with the case when $c = 1$, $\gamma$ decreases with increasing $L$ and $m$ when the other parameter is fixed, for any value of $c$. Further, when $\kappa$ is fixed and $m$ or $L$ are made twice as large, $\gamma$ decreases by a factor of two, for any value of $c$. This can be seen by comparing Fig.~\ref{fig:gd-varying-na}(a) with Fig.~\ref{fig:gd-varying-na}(b), and comparing Fig.~\ref{fig:gd-varying-na}(c) with Fig.~\ref{fig:gd-varying-na}(d). Thus, $\gamma$ varies inversely with $m$ and $L$ for a fixed $\kappa$, and continues to depend separately on these parameters. The value of $\gamma$ increases as $c$ increases, i.e. the algorithm is more sensitive to noise when the interval $\mathcal{A}$ is larger. These results are also consistent with the bound derived in Proposition \ref{prop:NA} for different values of $c$.

A few remarks are in order. In both Fig.~\ref{fig:gd-varying-iter-conv} and Fig.~\ref{fig:gd-varying-na}, the guarantees for $c > 1$ are generally worse than the guarantees for $c = 1$, a constant step size. The results show that convergence is slower, and the algorithm is more sensitive to noise. This is primarily a consequence of the setup discussed in Section \ref{sec:Varying}. The approach {we discuss in Section \ref{sec:Implementation}} is inherently conservative, and the results demonstrate worse guarantees than the constant step size case, especially when the step size is allowed to vary over a larger set. However, in practice, line search algorithms generally perform better (and not worse as predicted by our analysis) than their constant step size counterparts in terms of convergence or noise amplification. While the results do not accurately represent this, future work must focus on improving the methods discussed in Section \ref{sec:Varying}. In particular, it is worth exploring how line search algorithms (such as a Wolfe line search) can be more accurately characterized so that the IQC approach may lead to improved guarantees. A difficulty with this, as mentioned in Section \ref{sec:Varying} is that the derivable range on the step size using a Wolfe line search can be much more conservative than the range used here, resulting in worse guarantees.

\section{Concluding Remarks} \label{sec:Conclusion}
In this article, we presented an analysis of gradient descent with varying step sizes, in terms of its convergence rate and noise amplification. Assuming a line search produces a step size in a given interval, the algorithm is modeled as an LPV system. 
{Using a technique for polytopic LPV systems and building on prior work in the IQC framework, we construct SDPs to certify convergence rates and the steady-state variance in iterate error.} 
Our condition provides a bound on the convergence rate when the step size is within a restricted set around $\alpha = 1/L$. Moreover, our condition recovers the corresponding gradient rate when the interval is a single point $1/L$. Further, the noise amplification metric depends on both parameters $m$ and $L$ individually, and varies inversely with them for a fixed condition number $\kappa$.

It is worth reiterating that the guarantees for a line search, i.e. for $c > 1$ are generally worse than those for a constant step size. This is not necessarily reflective of the practical performance of line search algorithms, and it is worth exploring how these algorithms can be more accurately characterized to obtain less conservative results. Another avenue for further work in this area is the use of dynamic IQCs or multiple IQCs that may reduce conservatism in the results. Numerous dynamic IQCs for convex functions (where $\Psi$ in Section \ref{sec:IQCs} is not simply a static map, but a dynamical system) are described in \cite{LL2016}, and the use of multiple such IQCs may reduce conservatism in both the convergence rate bound and the noise amplification metric. 
Further avenues for future work include the analysis of time-varying step sizes in accelerated and stochastic gradient algorithms.

\appendix
\section{Proof of Analytical Results}
\subsection{Proof of Proposition \ref{prop:CR}} First, note that \eqref{eq:VaryingConvLMI} is a family of two LMIs at $\underline{\alpha}$ and $\overline{\alpha}$, based on the results in \eqref{eq:LMISimplified-1}, \eqref{eq:LMISimplified-2} and \eqref{eq:LMISimplified-3}. For a given $\alpha$, we know from the constant step size case that the convergence rate satisfies \eqref{eq:const-rho}.
Let $\underline{\rho}$ and $\overline{\rho}$ be the convergence rate corresponding to the two step sizes $\underline{\alpha}$ and $\overline{\alpha}$. Then, the convergence rate from \eqref{eq:VaryingConvLMI} can be written as:
\begin{equation}
    \rho_{\mathrm{varying}} = \max\{\underline{\rho}, \overline{\rho}\}
\end{equation}
since this is the smallest step size for which \eqref{eq:VaryingConvLMI} is feasible at both $\underline{\alpha}$ and $\overline{\alpha}$, {with a common solution $P = 1$ which is discussed in \cite{LL2016}.} Since $\underline{\alpha} = \frac{1}{cL} \leq \frac{2}{L+m}$, using \eqref{eq:const-rho}:
\begin{equation} \label{eq:under-rho}
    \underline{\rho} = 1 - \frac{m}{cL} = 1 - \frac{1}{c\kappa}.
\end{equation}
For $\overline{\rho}$, note that $\overline{\alpha} = \frac{c}{L} \geq \frac{2}{L+m}$ if $\kappa \leq \frac{c}{2-c}$ and $c < 2$. Then, using \eqref{eq:const-rho} with $c < 2$:
\begin{equation} \label{eq:over-rho}
    \overline{\rho} =
    \begin{dcases}
        c-1, \hspace{1em} &\kappa \leq \frac{c}{2-c}, \\
        1 - \frac{c}{\kappa}, &\kappa > \frac{c}{2-c}
    \end{dcases}.
\end{equation}
Comparing \eqref{eq:under-rho} and \eqref{eq:over-rho}, first note that $\frac{1}{c(2-c)} \leq \frac{c}{2-c}$ for $1 \leq c < 2$. When $\kappa \leq \frac{1}{c(2-c)} \leq \frac{c}{2-c}$, $\rho_{\mathrm{varying}} = \max\left\{1-\frac{1}{c\kappa}, c-1\right\} = c-1$ in this range of $\kappa$. Next, when $\frac{1}{c(2-c)} \leq \kappa \leq \frac{c}{2-c}$, $\rho_{\mathrm{varying}} = \max\left\{1-\frac{1}{c\kappa}, c-1\right\} = 1-\frac{1}{c\kappa}$ in this range of $\kappa$. Finally, when $\kappa > \frac{c}{2-c}$, $\rho_{\mathrm{varying}} = \max\left\{1-\frac{1}{c\kappa}, 1-\frac{c}{\kappa}\right\} = 1-\frac{1}{c\kappa}$ for $c \geq 1$. Using the above facts, we have:
\begin{equation} \label{eq:app-final-rho}
    \rho_{\mathrm{varying}} = 
    \begin{dcases}
        c-1, \hspace{1em} &\kappa \leq \frac{1}{c(2-c)} \\
        1-\frac{1}{c\kappa}, &\kappa > \frac{1}{c(2-c)} 
    \end{dcases}
\end{equation}
for $1 \leq c < 2$, which completes the proof of Proposition \ref{prop:CR}. \hfill \qedsymbol

\subsection{Proof of Proposition \ref{prop:NA}} {From \eqref{eq:SDP-gamma}, we require $\gamma_{\mathrm{varying}}^{2} \geq \mathrm{tr}\left(B_{w}^{T}(\alpha)PB_{w}(\alpha)\right)$ at $\underline{\alpha}$ and $\overline{\alpha}$, where $P$ is a common solution of the LMI \eqref{eq:VaryingNALMI} at these two points. Using the dimensionality reduction argument described in \cite{LL2016}, we can restrict our search to those $P$ of the form $P = P_0\otimes I_n$, where $P_0$ is $1\times 1$. Substituting $P$ and $B_w(\alpha)$, we require $\gamma_{\mathrm{varying}}^{2} \geq \alpha^2nP_0$ at both $\underline{\alpha}$ and $\overline{\alpha}$, or 
\begin{equation} \label{eq:gamma-P-n}
	\gamma_{\mathrm{varying}}^{2} \geq \overline{\alpha}^2nP_0
\end{equation}
since $\overline{\alpha} \geq \underline{\alpha}$. Let $\underline{P}_0$ and $\overline{P}_0$ be the two solutions at the points $\underline{\alpha}$ and $\overline{\alpha}$ respectively. Using \eqref{eq:gamma-P-n}, we then require:
\begin{equation} \label{eq:gamma-pstar}
\gamma_{\mathrm{varying}}^{2} \geq \left(\gamma^{*}_{\mathrm{varying}}\right)^2 = \overline{\alpha}^2nP_{0}^{*},~~~P_{0}^{*} = \max\{\underline{P}_0, \overline{P}_0\}.
\end{equation}
For a given $\alpha$, the solution $P_0$ to \eqref{eq:VaryingNALMI} is given by:
\begin{equation} \label{eq:Sol_P}
P_0 = \begin{dcases} \frac{1}{2\alpha m-\alpha^2m^2}, &\alpha \leq \frac{2}{L+m}, \\
\frac{1}{2\alpha L - \alpha^2L^2}, &\alpha \geq \frac{2}{L+m} \end{dcases},
\end{equation}
which follows from similar arguments to \cite{MRJ19}. Substituting $\underline{\alpha} = \frac{1}{cL}$ and $\overline{\alpha} = \frac{c}{L}$, 
\begin{align} 
\underline{P}_0 &= \frac{1}{\frac{2}{c\kappa}-\frac{1}{(c\kappa)^2}}, \label{eq:under-P} \\
\text{and }~\overline{P}_0 &= \begin{dcases} \frac{1}{2c-c^2}, &\kappa \leq \frac{c}{2-c}, \\
\frac{1}{\frac{2c}{\kappa} - \left(\frac{c}{\kappa}\right)^2}, &\kappa > \frac{c}{2-c} \end{dcases}. \label{eq:over-P}
\end{align}
This follows from noticing $\overline{\alpha} = \frac{c}{L} \geq \frac{2}{L+m}$ if $\kappa \leq \frac{c}{2-c}$ and $c < 2$. Note that $\underline{P}_0$ and $\overline{P}_0$ are both positive since $1 \leq c < 2$ and $\kappa \geq 1$. Finally, it is easily shown that $\overline{P}_0 \geq \underline{P}_0$ for all $\kappa$ and $c < 2$, and thus $P_{0}^{*} = \overline{P}_0$. Using \eqref{eq:gamma-pstar} and simplifying,
\begin{equation} \label{eq:app-final-gamma}
    \gamma_{\mathrm{varying}} \geq \gamma^{*}_{\mathrm{varying}} = 
    \begin{dcases}
        \frac{c}{L}\sqrt{\frac{n}{2c-c^2}}, \hspace{1em} &\kappa \leq \frac{c}{2-c} \\
        \frac{c}{m}\sqrt{\frac{n}{2c\kappa - c^2}}, &\kappa > \frac{c}{2-c}
    \end{dcases}
\end{equation}
for $1 \leq c < 2$, which completes the proof of Proposition \ref{prop:NA}.} \hfill \qedsymbol

{
\bibliographystyle{IEEEtran}
\bibliography{refs}
}

\end{document}